\documentclass[review]{elsarticle}
\journal{Journal of Computational Physics}

\usepackage{fullpage}
\usepackage[hang,flushmargin]{footmisc}
\usepackage{graphicx}
\usepackage{caption}
\usepackage{subcaption}
\usepackage{amssymb,amsmath}
\usepackage{bm}

\begin{document}

%%===============================
%% MACROS
%%===============================

\newcommand{\Thermochimica}{\textsc{thermochimica}}
\newcommand{\Optima}{\textsc{optima}}
\newcommand{\FactSage}{\textsc{factsage}}
\newcommand{\OptiSage}{\textsc{optisage}}
\newcommand{\FACT}{\textsc{fact}}
\newcommand{\Pandat}{\textsc{pandat}}
\newcommand{\ChemSage}{\textsc{chemsage}}
\newcommand{\ThermoCalc}{\textsc{thermocalc}}
\newcommand{\BINGSS}{\textsc{bingss}}
\newcommand{\Parrot}{\textsc{parrot}}
\newcommand{\FAST}{\textsc{fast}}
\newcommand{\MPMFAST}{\textsc{mpm-fast}}
\newcommand{\comsol}{\textsc{comsol}}

%%===============================
%% FRONT MATTER
%%===============================

\begin{frontmatter}

\title{A Jacobian Free Deterministic Method for Solving Inverse Problems}

\author[fesns]{M.H.A. Piro\corref{cor1}}

\ead{markus.piro@uoit.ca}
\cortext[cor1]{Corresponding author.}

\author[rmcc,cnlfd]{J.S. Bell}
\author[fesns]{M. Poschmann}
\author[cnlct]{A. Prudil}
\author[rmcc]{P. Chan}

\address[fesns]{Faculty of Energy Systems and Nuclear Sciences, Ontario Tech University, Oshawa, ON, Canada}
\address[rmcc]{Department of Chemistry and Chemical Engineering, Royal Military College of Canada, Kingston, ON, Canada}
\address[cnlfd]{Fuel Development Branch, Canadian Nuclear Laboratories, Chalk River, ON, Canada}
\address[cnlct]{Computational Techniques Branch, Canadian Nuclear Laboratories, Chalk River, ON, Canada}

%%------------------------------------------
%% ABSTRACT
%%------------------------------------------

\begin{abstract}
\label{KeyAbstract}
\noindent An effective numerical method is presented for optimizing model parameters
that can be applied to any type of system of non-linear equations and any number of data-points,
which does not require explicit formulation of the objective function or its partial derivatives.
The numerics are reduced to solving a non-linear least squares problem, which
uses the Levenberg-Marquardt algorithm and the Jacobian is approximated by applying
rank-one updates using Broyden's method. An advantage of this methodology over conventional
approaches is that the partial derivatives of the objective function do not have to be analytically
calculated. For instance, there may be situations where one cannot formulate the partial derivatives,
such as cases involving an objective function that itself contains a nested optimization problem.
Moreover, a line search algorithm
is also described that ensures that the Armijo conditions are satisfied and that convergence is assured,
which makes the success of the approach insensitive to the initial estimates of the model parameters.
 The foregoing numerical methods are described with respect to the development of the \Optima{}
software to solve inverse problems, which are reduced to non-linear least squares problems.
This computational approach
has proven to be particularly useful at solving inverse problems of very complex physical
models that cannot be optimized directly in a practical way.

\end{abstract}

\end{frontmatter}

%%===============================
%% INTRODUCTION
%%===============================

\section{Introduction}
\label{KeyIntro}

Frequently in science and engineering one is required to solve an inverse problem, which is the optimization of model parameters to yield model predictions consistent with experimental observations. Often, this results in solving a Non-Linear Least Squares (NLLS) problem, which requires the calculation of a set of partial derivatives of the objective function to be minimized. The Levenberg-Marquardt Algorithm (LMA) is a popular and effective method to solve this type of problem, which requires a linear system involving a Jacobian matrix to be solved. A great challenge in some inverse problems is that the partial derivatives of the objective function cannot be calculated directly in a practical way or it is undesirable -- for whatever reason -- to do so. For example, there may be situations where the NLLS problem requires solving a system of nested non-linear problems, or depends on stochastic outputs. Also, there may be situations where it is desirable to solve a NLLS problem coupled to a black-box code (e.g., closed-source commercial software). Under these circumstances, it may be desirable to solve a NLLS problem without calculating the partial derivatives explicitly in order to construct a Jacobian matrix.

The approach described herein solves the NLLS problem using LMA, but what is unique is that the Jacobian matrix is approximated using rank-one updates using Broyden's method. This eliminates the need to compute the partial derivatives directly making the approach quite versatile and ``Jacobian free'', since the Jacobian matrix is not computed explicitly but rather approximated by a Broyden matrix. Furthermore, a line search algorithm is described whereby the step-length is calculated in a manner that ensures the Armijo conditions are satisfied, which promotes numerical robustness. These numerical methods have been incorporated in the development of the software \Optima{}.

Two practical example problems are described that make use of \Optima{}. The first involves the development of a thermodynamic model with some commercial software in conjunction with \Optima{}, while the second describes the development of a fission gas diffusion model using a multi-physics code with \Optima{}. For both cases, descriptions are provided on the development of each model with respect to an appropriate set of experimental data. The numerical methods that broadly describe solving a non-linear least squares problem are presented in \S \ref{KeyLeastSquares} and the application of these methods is described in \S \ref{KeyImplementation}. Finally, \S \ref{KeyExample} describes the application of \Optima{} to solving two inverse problems of practical engineering interest.

%%===============================
%% NON-LINEAR LEAST SQUARES
%%===============================

\section{Non-Linear Least Squares Methods}
\label{KeyLeastSquares}

A common numerical problem is to fit a parameterized function to a set of data, which is accomplished by minimizing the sum of squares of the residuals represented by the functional value and the corresponding datum of individual terms. The method of least squares reduces to minimizing an objective function given by
\begin{equation}
	\label{EqObjFunc}
	S(\mathbf{x},\bm{\beta}) = \tfrac{1}{2} \sum_{i=1}^{m} \left(r_i(x_i,\bm{\beta})\right)^2
\end{equation}

\noindent where $r_i(x_i,\bm{\beta})$ is the residual of the $i$th term, given by
\begin{equation}
	\label{EqResidual}
	r_i(x_i,\bm{\beta}) = y_i - f(x_i, \bm{\beta})
\end{equation}

\noindent where $x_i$ and $y_i$ are the $m$ independent and dependent variables of the empirical datum, respectively, $\bm{\beta}$ is the vector of $n$ unknown parameters, $f$ is the functional model to be fitted, and $S$ is the objective function to be minimized. By definition of a least squares problem, the system is over-determined and $m \ge n$. For situations where $f$ is linear, a linear-least squares algorithm may be used, which can be solved directly. On the contrary, an iterative approach is required when $f$ is non-linear, which requires careful consideration and will be the focus of attention for the remainder of this section.

%% ------------------------------------------
%% Levenberg-Marquardt Algorithm
%% ------------------------------------------

\subsection{The Levenberg-Marquardt Algorithm}
\label{KeyLMA}

LMA is a popular method for solving non-linear least squares problems and it is employed in the current work. Unlike other methods, such as the Gauss-Newton method, LMA replaces a line search for a trust region strategy~\cite{Nocedal06}. Utilization of a trust region circumvents numerical issues when the Jacobian is rank deficient, which often plagues the Gauss-Newton method~\cite{Nocedal06}. One could interpret the LMA approach as an interpolation between the Gauss-Newton method and the method of gradient descent. As a result of this hybrid approach, the LMA is more robust than the Gauss-Newton method, even when the initial estimates are far from the final solution~\cite{Nocedal06}.

The LMA method reduces to solving the following system of linearized
equations~\cite{Nocedal06}:
\begin{equation}
	\label{EqLMA}
	(\mathbf{J^TJ} + \lambda \, \mathrm{diag} (\mathbf{J^TJ))p} = -\mathbf{J^T r}
\end{equation}

\noindent where $\mathbf{J} \in \Re^{m \times n}$ is the Jacobian matrix, $\mathbf{p} \in \Re^{n}$ is the direction vector, $\mathbf{r} \in \Re^{m}$ is the residual vector, $\mathbf{T}$ is the transpose operator, and $\lambda$ is a non-negative damping factor (scalar). The Jacobian contains the first-order partial derivatives of the objective function, $S$, with respect to the unknown variables. Applying a small value of $\lambda$ results in a Gauss-Newton update, while a relatively large value puts the system in gradient descent. In practice, it is useful to be in gradient descent (relatively large $\lambda$) when the system is very far from the expected minimum and to progressively move towards a Gauss-Newton step (relatively small $\lambda$) as the system begins to converge.

The coupled system of linearized equations can also include a weighting matrix, $\mathbf{W} \in \Re^{m \times m}$, which is a diagonal matrix containing weighting factors. These values may be used, for example, to represent errors or uncertainties in experimental measurements. Incorporating the weighting matrix $\mathbf{W}$ in equation (\ref{EqLMA}) yields
\begin{equation}
	\label{EqLMAW}
	(\mathbf{J^TWJ} + \lambda \, \mathrm{diag} (\mathbf{J^TWJ))p} = -\mathbf{J^T W r}
\end{equation}

The objective of LMA is to solve for the direction vector $\mathbf{p}$ in either equation (\ref{EqLMA}) or (\ref{EqLMAW}) in order to update the parameter vector $\bm{\beta}$. The system of linear equations containing the Jacobian matrix must be solved in order to compute the direction vector, $\mathbf{p}$; however, the calculation of $\mathbf{J}$ is not always straightforward and often requires special attention, which will be more carefully examined in the next section.

%% ------------------------------------------
%% Estimating the Jacobian
%% ------------------------------------------

\subsection{Estimating the Jacobian Matrix}
\label{KeyJacobian}

One may be able to analytically compute all of the partial derivatives required in order to construct the Jacobian matrix for a number of NLLS problems; however, this can be very challenging or impractical for some types of numerical problems that involve very complex functional models. The current approach does not rely on the ability to analytically calculate $\mathbf{J}$. Instead, Broyden's method is utilized to estimate the Jacobian by applying rank-one updates to an approximation of $\mathbf{J}$, which will heretofore be referred to as the Broyden matrix, $\mathbf{B}$. Thus, a linear model of the objective function is created that resembles the exact Jacobian, whereby $\mathbf{B} \cong \mathbf{J}$. The update to $\mathbf{B}$ at iteration $k$ is defined by the observed change in the functional vector, $\mathbf{t}_k = \mathbf{r}_k - \mathbf{r}_{k-1}$, with respect to changes in the system variables, $\mathbf{s}_k = \bm{\beta}_k - \bm{\beta}_{k-1}$, given by~\cite{Nocedal06}
\begin{equation}
	\label{EqBroyden}
	\mathbf{B}_k = \mathbf{B}_{k-1} + \frac{(\mathbf{t}_k - \mathbf{B}_{k-1}
		\mathbf{s}_k ) \mathbf{s_k^T}} {\| \mathbf{s}_k \| ^2}
\end{equation}

Once $\mathbf{B}_k$ has been updated, the direction vector $\mathbf{p}$ can be solved from either equation (\ref{EqLMA}) or (\ref{EqLMAW}) using a linear equation solver and the system variables can be updated with an appropriate step length, $\alpha_k \in (0,1]$. The update to the system parameters is given by
\begin{equation}
	\label{EqUpdate}
	\bm{\beta}_{k} = \bm{\beta}_{k-1} + \alpha_k \mathbf{p}_k
\end{equation}

As is frequently the case in many numerical optimization problems, the success of the approach rests on an effective determination of a step length. The step length is selected to ensure a sufficient decrease of the functional norm of equation (\ref{EqObjFunc}) by enforcing the Armijo conditions~\cite{Nocedal06}, given by
\begin{equation}
	\label{EqArmijo}
	\| \mathbf{r}( \mathbf{x}, \bm{\beta}_{k-1}
		+ \alpha_k \mathbf{p}_k) \|_2 \le \| \mathbf{r}( \mathbf{x},\bm{\beta}_{k-1} ) \|_2
		+ c \alpha_k \mathbf{r^T}(\mathbf{x}, \bm{\beta}_{k-1})\mathbf{p}_k
\end{equation}

\noindent where $c \in (0,1)$ and is often taken to be arbitrarily small in comparison to the functional norm. In the event that the Armijo conditions are not satisfied, the steplength is progressively halved and the system is reiterated until the inequality given by equation (\ref{EqArmijo}) is satisfied or a minimum step length is reached (e.g., $10^{-4}$). Furthermore, it is also common practice to constrain $\alpha$ to prohibit the system from departing from the feasible region. This can be accomplished, for example, by computing $\alpha$ such that the change of an independent variable is half-way between the domain of the feasible region and the previously estimated value.

Convergence can be achieved when a specific set of criteria is met. For instance, one criterion is when the relative change of all parameters is below a pre-specified limit, defined by
\begin{equation}
	\label{EqConvergence}
	\mathrm{max}(| p_j / \beta_j|) < \epsilon
\end{equation}

\noindent where $\epsilon$ is a pre-defined numerical tolerance that can be controlled by the user (e.g., $10^{-3}$). Similarly, one could use the magnitude of $\mathbf{p}$, a maximum number of iterations, and perhaps other conditions.

%%===============================
%% NUMERICAL IMPLEMENTATION
%%===============================

\section{Numerical Implementation}
\label{KeyImplementation}

The previous section described a general method to solve a non-linear least squares problem that does not require the analytical calculation of the partial derivatives to construct the Jacobian matrix. This section provides details on the numerical implementation of the foregoing method to optimize
model parameters.  All that is required to optimize a set of parameters -- regardless of what the parameters represent -- is the calculation of the residual vector, $\mathbf{r}_k$, corresponding to a particular set of parameters, $\bm{\beta}_k$, at iteration $k$.

As is typically the case in quasi-Newton methods, two sets of initial values are required to initialize the computation. Initial values of the parameter vector $\bm{\beta}_0$ may be provided by the user to accelerate the process or are otherwise taken to be equal to zero. A small perturbation from the initial set of $\bm{\beta}_0$ is applied for the second estimate to initialize the solution and the iteration procedure commences. The method will then determine a unique combination of parameters that minimizes the sum of squares of the residual terms for a given set of variables.

The following steps summarize the procedure required to perform model parameter optimization:

\begin{enumerate}
	\item Initialize $\bm{\beta}_0$ (e.g., default value of zero) and compute $\mathbf{r}_0$
		for all datum specified.
	\item Perturb $\bm{\beta}_0$ by a small quantity to produce $\bm{\beta}_1$ and
		compute $\mathbf{r}_1$.
	\item Compute $\mathbf{t}_k$ and $\mathbf{s}_k$ using the previous two sets of estimates,
		and update $\mathbf{B}_k$ using equation (\ref{EqBroyden}).
	\item Compute $\mathbf{p}_k$ from equation (\ref{EqLMA}) or equation (\ref{EqLMAW})
		(depending on whether weighting factors are used) using a linear equation solver.
	\item Update $\bm{\beta}_{k}$ using equation (\ref{EqUpdate}) with a steplength $\alpha_k$
		that ensures that the Armijo conditions are satisfied (i.e., equation (\ref{EqArmijo})).
	\item Check for convergence (i.e., equation (\ref{EqConvergence})). If converged, stop;
		if not, return to \#3 and reiterate.
\end{enumerate}

The first two steps are solely used to initialize the calculation, while steps 3-8 constitute an iterative cycle. The perturbation applied in \#2 may be problem specific, but a small perturbation of $\approx1\%$ has been found to work well for most problems.

%%===============================
%% EXAMPLES
%%===============================

\section{Example Problems}
\label{KeyExample}

Two example problems are presented that use the aforementioned methodology to develop material models. The first example describes the development of a thermodynamic model of hypostoichiometric thoria, which is a pertinent system to the nuclear community as a candidate fuel material (or component in a multi-component material) for future potential industrial use. The second example describes the development of a fission gas diffusion model for use in a nuclear fuel performance model. Fortuitously, both examples involve thoria, whereby the latter considers thoria-plutonia as a mixed oxide.

The examples involve the computation of some non-linear problem using two different commercial software packages. For both cases, it was desirable to solve the inverse problem for some material behaviour that could not be solved directly by the respective software. Since both software packages are closed-source, the typical user does not have access to the source code and cannot make any modifications if one wishes to perform some optimization process. In these situations, there is a need to be able to solve a NLLS problem when one cannot solve the partial derivatives, since one cannot access them directly through the commercial software. Furthermore, for both cases, the minimization of $S$ requires a nested calculation of a non-linear system of equations.

%% ------------------------------------------
%% Thermodynamic model development
%% ------------------------------------------

\subsection{Example I: Development of a Thermodynamic Model of Hypostoichiometric Thoria}

The development of thermodynamic models rests on the optimization of function parameters to yield consistent predictions with experimental measurements and has been formalized by the CALPHAD method~\cite{Lukas07}. Specific function parameters include molar heat capacity and excess molar Gibbs energy of mixing between constituents in a non-ideal solution phase. The optimization of these parameters are typically made with respect to experimental measurements of the molar heat capacity, thermochemical activities of various components, and phase boundary data. A coefficient of the residual vector may represent, for example, the difference between measured and predicted values of the natural logarithm of the thermochemical activity of a particular chemical species. One could use the resulting models to perform a number of predictions, such as phase equilibria calculations to construct a phase diagram, which are useful in a number of areas of interest in materials science and many engineering disciplines.

Several commercial codes are available to assist the user in developing thermodynamic models -- including \FactSage{}~\cite{Bale02}, \Pandat{}~\cite{Cao09}, \ThermoCalc{}~\cite{Sundman85}, and \BINGSS{}~\cite{Lukas77} -- that have been successful in expediting the model development process. The success of these codes is largely due to their versatility and convenience, which circumvents the painstaking process of developing a model by manipulating model parameters by trial-and-error or by performing optimization calculations by hand. However, the success of these approaches often depends on the user's initial estimate of the assemblage of coefficients to be optimized.

The current work builds on the knowledgebase put forth by existing thermodynamic optimization methods by incorporating more sophisticated numerical methods that are commonly employed in other applications of applied mathematics. Thus, the intent of the current work is to provide a numerical approach that is versatile, convenient and robust, whereby the success of the calculation is independent of the initial estimates. Furthermore, previously published articles related to this topic have focused on the application of computational tools to develop thermodynamic databases, whereas the current article will focus on the numerical methods themselves with greater mathematical rigor.

The objective of optimizing thermodynamic model parameters reduces to solving a NLLS problem, which can be solved by a variety of numerical approaches, each of which have their own unique characteristics. Existing thermodynamic model optimizers employ different numerical algorithms to solve the NLLS problem. For instance, \OptiSage{} (i.e., \FactSage{}) and its predecessor \ChemSage{} have enjoyed success using the Bayesian approach~\cite{Konigsberger91,Konigsberger95}, while \FACT{}~\cite{Bale83}, \Parrot{} (i.e., \ThermoCalc{}), and \BINGSS{} have relied on the Gauss-Newton approach~\cite{Lukas07,Lukas77}. In these implementations, the partial derivatives of the objective function appear to be approximated to construct the Jacobian matrix; however, very little information has been disclosed in the open literature, presumably because they have been largely developed for commercial software and constitute trade secrets.

For the case of thermodynamic model parameter optimization, the objective function depends on the predicted assemblage of stable phases, the types of thermodynamic models that are employed (e.g., regular substitutional model or compound energy formalism), the types of parameters that are altered (e.g., coefficients of the standard reference Gibbs energy terms or excess Gibbs energy of mixing terms), and the empirical data to be optimized. Furthermore, the complexity of these calculations increases with higher order systems.

The difficulty in optimizing thermodynamic model parameters is that it requires the calculation of a nested subset of optimization problems. In other words, each of the $m$ terms involved in minimizing $S$ from equation (\ref{EqObjFunc}) requires the minimization of the integral Gibbs energy of a closed isothermal-isobaric system. The Lagrangian function to be minimized within the nested subset of calculations takes the form~\cite{Eriksson73,Jansson84}:
\begin{equation}
	\label{EqGEM}
	L = \sum_{\phi=1} n_{\phi} g_{\phi} - \sum_{j=1} \Gamma_j \left( b_j - \sum_{\phi=1} n_{\phi} \nu_{\phi,j} \right)
\end{equation}

\noindent where $n_{\phi}$ and $g_{\phi}$ represent the number of moles and the molar Gibbs energy of phase $\phi$, respectively, $\Gamma_j$ is the Lagrangian multiplier corresponding to system component $j$ (equivalently, the chemical potential of $j$), and$\nu_{\phi,j}$ is the effective stoichiometry coefficient of component $j$ in phase $\phi$. The minimization of $L$ is subject to linear and nonlinear equality and inequality constraints, which is handled in a separate manner from $S$. Furthermore, the $g_{\phi}$ function is non-linear and often non-convex, which requires a global optimization algorithm to ensure convergence~\cite{Piro16} and the active set of constraints often have to be updated throughout the iterative cycle~\cite{Piro17}. It is for these reasons that one cannot compute the partial derivatives of equation (\ref{EqObjFunc}) involving $L$ analytically.

Once $L$ has been minimized, one can use the results to formulate equation (\ref{EqObjFunc}). For example, one can compute the residuals of the chemical potentials of a particular component that are measured and predicted (i.e., $\Gamma_j$). In this work, model parameters of the ThO$_{2-x}$ phase were optimized with respect to a set of experimentally measured oxidation data. The specific model is based on the Compound Energy Formalism (CEF)~\cite{Hillert01} with three sublattices with the intention of being fully compatible with the model developed by Gu{\'{e}}neau \textit{et al.}~\cite{Gueneau11} and the Thermodynamics of Advanced Fuels International Database (TAF-ID)~\cite{Gueneau15}. In the CEF model, species are represented by a unique combination of a single constituent on each sublattice, which are often referred to as ``compound end members''.

The general equation for the molar Gibbs energy of solution phase $\lambda$, where $\lambda \subset \phi$, that is represented by CEF is given by:
\begin{equation}
	\label{EqCEF}
	g_{\lambda} = \sum_{i=1}^{N_{\lambda}} g_{i(\lambda)}^o \prod_{s=1}^{N_s} y_{i(s)}
		+ \mathrm{RT} \sum_{s=1}^{N_s} a_s \sum y_{i(s)} \mathrm{ln} (y_{i(s)} )
		+ g_{\lambda}^{ex}
\end{equation}

\noindent where $i$ is the species index, $N_{\lambda}$ is the number of species in solution phase $\lambda$, $g_{i(\lambda)}^o$ is the standard molar Gibbs energy of pure species $i$, $y_{i(s)}$ represents the site fraction of a particular constituent on sublattice $s$ corresponding to species $i$, $R$ is the ideal gas constant, $T$ is the absolute temperature, $\nu_s$ is the stoichiometric coefficient of sublattice $s$, and $g_{\lambda}^{ex}$ is the molar excess Gibbs energy of mixing.

The two variables in equation (\ref{EqCEF}) that are optimized are contained within $g_{i(\lambda)}^o$ and $g_{\lambda}^{ex}$. The former is given by the following equation:
\begin{equation}
	\label{EqSGE}
	g_{i(\lambda)}^o = h_{i(\lambda)}^o + \int_{T_o}^T c_p dT - T \left( s_{i(\lambda)}^o + \int_{T_o}^T \frac{c_p}{T} dT \right)
\end{equation}

\noindent where $h_{i(\lambda)}^o$ and $s_{i(\lambda)}^o$ are the standard enthalpy and entropy of formation at standard state, respectively, and $c_p$ is the temperature dependent molar heat capacity at constant pressure. The molar excess Gibbs energy is generally given as:
\begin{equation}
	\label{EqEG}
	g_{\lambda}^{ex} = \sum_{z=1} \sum_{n=0} {}^{n}L_z (y_a - y_b)^{n} \prod y_{i(\lambda)}^z
\end{equation}

\noindent where $z$ is a parameter index, $n$ is the order of the model parameter ${}^{n}L_z$ and $y_{i(\lambda)}^z$ represents the set of site fractions associated with model parameter $z$\footnote{To avoid confusion, it is worth acknowledging that the variable ``$L$'' is used to represent the Lagrangian function in equation (\ref{EqGEM}) and the non-ideal mixing parameter in equation (\ref{EqEG}). The notation used here simply follows convention, but can be distinguished by the use of super- and sub-scripts in the latter case.}. The variables $y_a$ and $y_b$ are the site fractions of constituents that are mixing for parameter $z$.

The parameters that are often optimized are $h_{i(\lambda)}^o$, $s_{i(\lambda)}^o$, $c_P$, and ${}^{n}L_z $. Note that since the compound end members are by definition related, the optimization of $h_{i(\lambda)}^o$, $s_{i(\lambda)}^o$, and $c_P$ necessarily affect multiple compound end members. This is often handled through the use of functions that directly affect the optimization of multiple model parameters. As an example, Gu{\'{e}}neau \textit{et al.} define all 12 compound end members in their model of UO$_{2\pm \textrm{x}}$ by functions of UO$_{1.5}$, UO$_2$, and UO$_{2.5}$ compounds (whereby the first and last are fictive)~\cite{Gueneau11}.

\begin{figure}[h]
	\includegraphics[width=12cm]{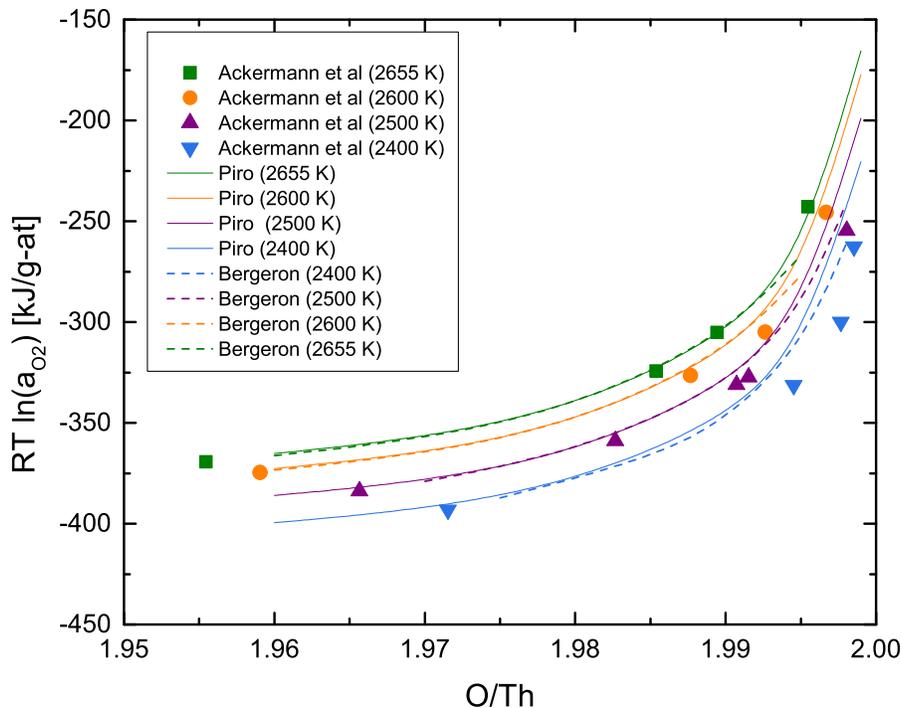}
	\centering
	\caption{Thermodynamic model predictions of oxygen partial pressure in equilibrium with ThO$_{2-x}$ are compared to experimental measurements by Ackermann \textit{et al.}~\cite{Ackermann80}}
	\label{FigThO2_O2}
\end{figure}

The commercial software \FactSage{} was used in performing thermodynamic calculations; however, the \OptiSage{} module within \FactSage{}, which is generally quite powerful, is unable to employ functions that apply to multiple model parameters simultaneously when applied to the compound end members. Instead, \Optima{} was used in conjunction with \FactSage{} to develop a model that did not invoke \OptiSage{}. The parameters were adjusted as input to \FactSage{} and the output was used as input to \Optima{}. The resulting model predictions of the thermochemical activity of oxygen are illustrated in Figure \ref{FigThO2_O2}, which compare very well to experimental data. Another model produced by Bergeron \textit{et al.}~\cite{BERGERON2018324} is compared in Figure \ref{FigThO2_O2}, which used \OptiSage{}; however, it required a greater number of tunable parameters. The purpose of this comparison is to demonstrate that \Optima{} works and that it can successfully be coupled to a blackbox code. To be clear, the problem is not as simple as fitting a polynomial to the experimental data points in Figure \ref{FigThO2_O2}, but rather optimizing model parameters that are used in numerous minimization calculations in equation (\ref{EqGEM}) that then yields the curves given in Figure \ref{FigThO2_O2}.

%% ------------------------------------------
%% Fission gas release
%% ------------------------------------------

\subsection{Example II: Development of a Fission Gas Diffusion Model in Nuclear Fuel}

The second example problem presented calculated the effective fission gas diffusion coefficient in thoria-plutonia nuclear fuel. The context of this work was done in reference to developing a nuclear fuel performance model of thoria-plutonia~\cite{Bell17} in order to expand the capabilities of the \FAST{} nuclear fuel performance computational framework~\cite{Prudil15a,Prudil15b}, which has been constructed within the \comsol{} multi-physics commercial software. Fuel performance codes calculate parameters that assess the fuel's fitness for service during irradiation.

An important and particularly challenging phenomenon to predict is the release of gaseous fission products (mainly atoms of Xenon and Krypton) from the fuel, which bear great importance to performance and safety. These Fission Gasses (FG) have low solubility in the fuel matrix and diffuse to the grain boundaries at a rate that is dependent on the local temperature, fission rate, grain-size, and fuel burnup\footnote{The term ``burnup'' is a measure of the amount of thermal energy released due to nuclear fission.}. This process is much better understood for UO$_2$ fuel, compared to ThO$_2$ and PuO$_2$ mixtures since UO$_2$ is much more widely used and studied. As a result, work on developing models for thoria based fuel performance models have started with a pre-existing urania model, and modifying the thermophysical and irradiation behaviour of the fuel as a surrogate~\cite{Long02,Lee04,Bjork15,Boer16}.

In the case of the FG diffusion model discussed here, the foundational research for UO$_2$ is the Booth model~\cite{Booth57}. In this formulation, the fuel is approximated as a collection of spherical grains with an effective radius in which the FG concentration, $C$ (at. m$^{-3}$), is calculated as a function of radial position according to the following partial differential equation~\cite{White83}:
\begin{equation}
	\label{EqFGDiffusion}
	\frac{\partial C}{\partial t} = \frac{D_0 b'}{b' + g_a} \nabla^2 C + P_{fg}
\end{equation}

\noindent where $t$ is time [s], $D_0$ is the single gas atom diffusion coefficient [m$^2$ s$^{-1}$], $b'$ is the intragranular resolution rate [s$^{-1}$], $g_a$ is the trapping rate [s$^{-1}$], and $P_{fg}$ is the volumetric FG production rate [at. m$^{-3}$ s$^{-1}$]. A detailed description is provided by Bell~\cite{Bell17}. The single gas atom diffusion coefficient can be further decomposed into individual mechanisms by:
\begin{equation}
	\label{EqFGDiffCoeff}
	D_0 = D_{thrm} + D_{irr} +D_{athrm}
\end{equation}

\noindent where $D_{thrm}$ is the diffusion coefficient due to thermally activated processes (i.e., $T > $ 1700 K), $D_{irr}$ is the diffusion coefficient due to irradiation induced vacancies (i.e., 1100 K $\leq T \leq 1700$ K), and $D_{athrm}$ is the diffusion coefficient due to athermal effects (i.e., $T <$ 1100 K).

To gain insight into the FG diffusion in grains of thoria-plutonia, it was assumed that the physical processes by which FG atoms are transported are the same as urania, although at different rates. This is accomplished by introducing weighting coefficients ($\mathbf{w} = w_1, w_2, w_3$) applied to the components of $D_0$ to account for the relative prevalence of each phenomena in thoria-plutonia compared to the urania baseline. This weighting is given by:
\begin{equation}
	\label{EqFGDiffCoeffBase}
	D_{0,Th} = w_1 D_{thrm} + w_2 D_{irr} + w_3 D_{athrm}
\end{equation}

In order to determine the weighting factors without direct experimental measurements of Fission Gas Release (FGR) rates through sweep-gas or in situ pressure measurements, the predicted cumulative end-of-life fission gas volumes as a function of the weightings were compared to post-irradiation examination measurements. The optimal weighting coefficients are those that allow \MPMFAST{} to best predict the experimental FGR measurements. This is an inverse problem to estimate the weightings (i.e., material properties) to obtain the best fit to experimentally measured FGR data in the absence of diffusivity data.

The weighting factors given by the vector $\mathbf{w}$ were computed by minimizing Equation (\ref{EqObjFunc}) using \Optima{}. Conventional Jacobian based optimization techniques cannot be applied since \MPMFAST{} cannot be analytically differentiated and is computationally expensive to evaluate precluding direct numerical differentiation to approximate the Jacobian.

Table \ref{TableFGR} shows the predicted FGR results after obtaining the optimized weightings with \Optima{} compared to the experimentally measured FGR data as a percentage of fission gas produced, along with the maximum calculated fuel temperature during irradiation. Experiments were performed at the Canadian Nuclear Laboratories (formerly Atomic Energy of Canada Limited) with the identifier BDL-422~\cite{Karam08,Karam10,Floyd11,Floyd16}. Measured FGR from the BDL-422 experiment was considerably lower than urania fuel under similar irradiation conditions, except for two cases that experience extended irradiation times and high temperature conditions.

\begin{table}[h]
\begin{center}
\caption{Simulated FGR results compared to experimental data from thoria-plutonia experiment BDL-422~\cite{Bell17}. }

\label{TableFGR}

  \begin{tabular}{|c|c|c|c}
	\hline
	\textbf{Case}			&\textbf{Predicted FGR [\%]}	&	\textbf{Measured FGR [\%]} 			\\
	\hline
	A	&	1.1 &		1.0	\\
	\hline
	B	&	0.0 &		0.3	\\
	\hline
	C	&	0.4 &		0.5	\\
	\hline
	D	&	0.0 &		0.1	\\
	\hline
	E 	&	0.0 &		0.2	\\
	\hline
  \end{tabular}
\end{center}
\end{table}

This datum was used to fit the weighting coefficients for $D_{athrm}$ and $D_{irr}$, which are the dominant diffusion processes below 1700 K, based on their predicted maximum temperature calculations. These simulation results are in reasonable agreement with the measured results from BDL-422 based on the uncertainty within the temperature calculation associated with the thoria-plutonia model~\cite{Bell17}. The second half of the dataset that experienced centreline temperatures greater than 1700 K from BDL-422 could be used to fit the weighting coefficient for $D_{thrm}$. While the possibility of modelling FG diffusion at elevated temperatures was examined by Bell~\cite{Bell17}, any calculation occurring above 1600 K carries an uncertainty associated with it.

The study and development of irradiation behaviour models of thoria-plutonia is still relatively new and has not yet fully matured~\cite{Bell17,Bjork15,Boer16} compared to those of urania and to a lesser extent urania-thoria~\cite{Long02,Lee04,Long04}. In the case of thoria-plutonia fuel, where experimental data is comparatively scarce, the inverse simulation approach enables the use of integral post-irradiation measurements to produce a preliminary model. With the varied irradiation conditions for each fuel element and the highly coupled nature of the fuel performance calculation, the computational expense of each simulation was significant. Each iteration of the model required approximately 2-3 h of simulation time to complete the set of cases from BDL-422. \Optima{} enabled the limited data set to be used effectively in the initial stages of determining the FG diffusion behaviour in thoria-plutonia, while greatly reducing the number of model iterations required to converge on the solution.

%%===============================
%% DISCUSSION
%%===============================

\section{Discussion}
\label{KeyDiscussion}

The numerical procedure presented in section \S \ref{KeyLeastSquares} describes a general method to optimize non-linear model parameters in \Optima{}. As discussed in \S \ref{KeyLeastSquares}, the current approach does not require the analytical derivation of the Jacobian, which is approximated by Broyden's method. In other words, the programming does not require specific formulation of the partial derivatives or a method of parsing equations. Thus, the programming is made much simpler by employing Broyden's method in comparison to analytically calculating the partial derivatives to construct the Jacobian matrix. Furthermore, the examples given in the previous section give evidence of the utility of this method to work with black box software -- specifically, \Optima{} was indirectly coupled to the commercial software \FactSage{} and \comsol{}.

Also worth noting is that for both cases described in the foregoing section, \Optima{} was initialized with really highly inaccurate estimates of both the independent variables and the Jacobian matrix. Recall that the Jacobian matrix is initially approximated by a Broyden unit matrix, which is a crude approximation of the true Jacobian. Also, the very first initial values of the independent variables were taken as zero. This gives evidence of the robustness of the technique, which has great value from a practical perspective since user intervention is removed.

The popular Broyden-Fletcher-Goldfarb-Shanno (BFGS) method, or any of its other closely related derivatives, is not used in the current approach because the BFGS method yields a matrix update that is necessarily symmetric and square. Recall that the Jacobian matrix in LMA is almost always rectangular, which is necessarily non-symmetric. Broyden's method, on the other hand, not only permits non-symmetric updates to the Broyden matrix, but it can also accommodate non-square matrices. One could, however, apply the BFGS algorithm to update the Hessian matrix, which is necessarily symmetric and square. The Hessian matrix is defined for the LMA by
 \begin{equation}
	\label{EqHessian}
	\mathbf{H} =\mathbf{J^TJ} + \lambda \, \mathrm{diag} (\mathbf{J^TJ)}
\end{equation}

\noindent which is comprised of two matrices, representing contributions from Gauss-Newton and gradient descent. The disadvantage of applying a quasi-Newton approach to approximate the Hessian matrix (i.e., $\mathbf{B} \cong \mathbf{H}$) in comparison to the approach described in \S \ref{KeyLeastSquares}, is that one would no longer be able to control $\lambda$. Therefore, one would not be able to traverse between gradient descent and Gauss-Newton methods, which is the primary motivation for using LMA. Furthermore, this is the same reason why the Sherman-Morrison~\cite{Sherman50} formulation of Broyden's method is not exploited in the current work. Although this formulation requires far less computational expense, the performance gains with using BFGS do not warrant the inability to control $\lambda$.

Alternatively, one could approximate $\mathbf{J}$ by a finite difference approach; however, this requires a relatively large number of additional function evaluations equal to $m \times n$ per global iteration. The advantage of computing $\mathbf{J}$ using finite differences in comparison to using rank-one updates using Broyden's method is that one could achieve a very good approximation of the Jacobian with an appropriately selected perturbation. Transtrum and Sethna~\cite{Transtrum12} have proposed a hybrid approach that employs both methods to approximate $\mathbf{J}$, whereby the Jacobian/Broyden matrix is periodically re-established using finite differences. Thus, computational performance and robustness may be improved for some problems.

In addition to the manner that the Jacobian is constructed, the current approach differs from traditional methods with the use of a backstepping line search algorithm. The use of this method to update the system variables ensures that the Armijo conditions are satisfied, which assures that the numerical system makes sufficient progress at each iteration and that it does not diverge. Thus, the method is particularly valuable when the iterative procedure begins with a relatively poor estimate of the model parameters, which eliminates the need for the user to provide good initial estimates.

Several existing thermodynamic tools -- such as \FactSage{}~\cite{Bale02}, \Pandat{}~\cite{Cao09}, \ThermoCalc{}~\cite{Sundman85} and \BINGSS{}~\cite{Lukas77} -- have enjoyed great success in numerical optimization for the purposes of model development. The two primary differences between the approach taken here and that of traditional methods is in the way that partial derivatives are handled and the use of a line search algorithm, which collectively promote numerical stability and flexibility. Future development plans include coupling of \Optima{} with \Thermochimica{}~\cite{Piro13} to
provide an effective open-source tool for thermodynamic model development.

In regards to fuel performance modelling, irradiated nuclear fuel is a complicated multidisciplinary system simultaneously involving aspects of solid mechanics, heat transfer, radiation transport, thermochemistry, and microstructural evolution. Experimentation with nuclear materials can be extremely expensive and time-consuming, requiring a great deal of analysis to maximize the value of each experiment. Inverse problems offer one way of extracting additional data from experimental data.

The codes used to simulate the fuel system are highly non-linear coupled multiphysics codes. In addition, many of these codes are legacy codes developed as stand-alone computer programs, which may include the physical models, discretization and solution methods into a single program. In the case of codes used in the safety analysis of nuclear power plants, there may be legal restrictions requiring the use of licensed, quality controlled codes that makes them impractical to modify for specific analyses. As a result, performing inverse studies to determine material properties may be difficult using traditional optimization techniques. Therefore, computationally efficient Jacobian free techniques, such as the subject of this paper, are value engineering tools.

The numerical approach presented in this work offers great flexibility in solving inverse problems. The user is ultimately responsible for making decisions pertinent to model development, who should also have a sufficient familiarity with the physical analogue to enable a reasonable representation of reality. A model development tool should never substitute scientific intuition, but rather to expedite the model development process.

%%===============================
%% CONCLUSION
%%===============================

\section{Conclusion}
\label{KeyConclusion}

A general numerical method is presented to optimize model parameters that is flexible and
robust. The exploitation of Broyden's
method in approximating the Jacobian matrix yields a method that is independent of the type of equation
being minimized, which does not require explicit formulation of the partial differentials.
Additionally, the utilization of the Levenberg-Marquardt algorithm together with an effective line search
algorithm provides a very robust approach for numerical optimization. By enhancing robustness through
these methods, the success of the overall approach becomes far less sensitive to the initial estimates.
Thus, the user is further removed from the optimization process, making the overall approach more
user friendly. The foregoing numerical methods have been incorporated into the \Optima{} software,
which can be used independently or directly coupled to other software as demonstrated through two
example problem. \\

%%===============================
%% ACKNOWLEDGEMENTS
%%===============================

\section{Acknowledgements}
\label{KeyAcknolwedgements}

J.S. Bell and P. Chan would like to thank several members of the Canadian Nuclear Laboratories (CNL) staff:
H. Hamilton (retired), J. Pencer, W. Richmond, I. Lusk (retired), A. Williams, and M. Floyd for their assistance and collaboration.
M.H.A. Piro thanks A. Bergeron (CNL), H.P. Gavin from Duke University, and M. Transtrum and J. Sethna from Cornell
University for helpful technical discussions. \\

The work on the thoria-plutonia model was funded under the Canada GEN-IV National Program. Funding to the Canada Gen-IV National Program was provided by Natural Resources Canada through the Office of Energy Research and Development, Atomic Energy of Canada Limited, and Natural Sciences and Engineering Research Council of Canada.
This research was undertaken, in part, thanks to funding from the Canada Research Chairs program of the Natural Sciences and Engineering Research Council of Canada.\\

%%===============================
%% REFERENCES
%%===============================

%\noindent \textbf{References}

\bibliographystyle{ieeetr} %% ieeetr.bst
\bibliography{PiroRef}		%% PiroRef.bib

\begin{thebibliography}{10}

\bibitem{Nocedal06}
J.~Nocedal and S.~Wright, {\em Numerical Optimization}.
\newblock Springer, 2006.

\bibitem{Lukas07}
H.~Lukas, S.~Fries, and B.~Sundman, {\em Computational Thermodynamics: The
  Calphad Method}.
\newblock Cambridge University Press, 2007.

\bibitem{Bale02}
C.~Bale, P.~Chartrand, S.~Degterov, G.~Eriksson, K.~Hack, R.~Mahfoud,
  J.~Melan{\c{c}}on, A.~Pelton, and S.~Peterson, ``{FactSage} thermochemical
  software and databases,'' {\em CALPHAD}, vol.~26, no.~2, 2002.

\bibitem{Cao09}
W.~Cao, S.-L. Chen, F.~Zhang, K.~Wu, Y.~Yang, Y.~Chang, R.~Schmid-Fetzer, and
  W.~Oates, ``{PANDAT} software with panengine, panoptimizer and
  panprecipitation for multi-component phase diagram calculation and materials
  property simulation,'' {\em CALPHAD}, vol.~33, no.~328-342, 2009.

\bibitem{Sundman85}
B.~Sundman, B.~Jansson, and J.-O. Andersson, ``The {T}hermo-{C}alc databank
  system,'' {\em CALPHAD}, vol.~9, no.~2, pp.~153--190, 1985.

\bibitem{Lukas77}
H.~Lukas, E.~Henig, and B.~Zimmermann, ``Optimization of phase diagrams by a
  least squares method using simultaneously different types of data,'' {\em
  CALPHAD}, vol.~1, no.~3, pp.~225--236, 1977.

\bibitem{Konigsberger91}
E.~K\"{o}nigsberger, ``Improvement of excess parameters from thermodynamic and
  phase diagram data by a sequential bayes algorithm,'' {\em CALPHAD}, vol.~1,
  no.~69-78, 15.

\bibitem{Konigsberger95}
E.~K\"{o}nigsberger and G.~Eriksson, ``A new optimization routine for
  {C}hem{S}age,'' {\em CALPHAD}, vol.~19, no.~2, pp.~207--214, 1995.

\bibitem{Bale83}
C.~Bale and A.~Pelton, ``Optimization of binary thermodynamic and phase diagram
  data,'' {\em Metallurgical Transactions B}, vol.~14B, pp.~77--83, 1983.

\bibitem{Eriksson73}
G.~Eriksson and E.~Rosen, ``General equations for the calculation of equilibria
  in multiphase systems,'' {\em Chemica Scripta}, vol.~4, pp.~193--194, 1973.

\bibitem{Jansson84}
B.~Jansson, {\em Computer Operated Methods for Equilibrium Calculations and
  Evaluation of Thermochemical Model Parameters}.
\newblock PhD thesis, Royal Institute of Technology, 1984.

\bibitem{Piro16}
M.~Piro and S.~Simunovic, ``Global optimization algorithms to compute
  thermodynamic equilibria in large complex systems with performance
  considerations,'' {\em Computational Materials Science}, vol.~118,
  pp.~87--96, 2016.

\bibitem{Piro17}
M.~Piro, ``Updating the estimated assemblage of stable phases in a gibbs energy
  minimizer,'' {\em CALPHAD}, vol.~58, pp.~115--121, 2017.

\bibitem{Hillert01}
M.~Hillert, ``The compound energy formalism,'' {\em Journal of Alloys and
  Compounds}, vol.~320, no.~2, pp.~161 -- 176, 2001.

\bibitem{Gueneau11}
C.~Gu{\'{e}}neau, N.~Dupin, B.~Sundman, C.~Martial, J.-C. Dumas,
  S.~Goss{\'{e}}, S.~Chatain, F.~D. Bruycker, D.~Manara, and R.~J. Konings,
  ``Thermodynamic modelling of advanced oxide and carbide nuclear fuels:
  Description of the {U-Pu-O-C} systems,'' {\em Journal of Nuclear Materials},
  vol.~419, no.~1-3, pp.~145--167, 2011.

\bibitem{Gueneau15}
C.~Gu{\'{e}}neau, S.~Goss{\'{e}}, A.~Quaini, N.~Dupin, B.~Sundman, M.~Kurata,
  T.~Besmann, P.~Turchi, J.~Dumas, E.~Corcoran, M.~Piro, T.~Ogata, R.~Hania,
  B.~Lee, R.~Kennedy, and S.~Massara, ``{FUELBASE}, {TAF-ID} databases and {OC}
  software: Advanced computational tools to perform thermodynamic calculations
  on nuclear fuel materials,'' in {\em The 7th European Review Meeting on
  Severe Accident Research}, 2015.

\bibitem{BERGERON2018324}
A.~Bergeron, D.~Manara, O.~Beneš, R.~Eloirdi, M.~Piro, and E.~Corcoran,
  ``Thermodynamic modelling of thoria-urania and thoria-plutonia fuels:
  Description of the th-u-pu-o quaternary system,'' {\em Journal of Nuclear
  Materials}, vol.~512, pp.~324--348, 2018.

\bibitem{Ackermann80}
R.~Ackermann and M.~Tetenbaum, ``High temperature thermodynamic properties of
  the thorium-oxygen system,'' {\em High Temperature Science}, vol.~13,
  pp.~91--105, 1980.

\bibitem{Bell17}
J.~Bell, {\em Thorium-based nuclear fuel performance modelling with
  {Multi-Pellet Material Fuel and Sheath Modelling Tool}}.
\newblock PhD thesis, Royal Military College of Canada, 2017.

\bibitem{Prudil15a}
A.~Prudil, B.~Lewis, P.~Chan, and J.~Baschuk, ``Development and testing of the
  {FAST} fuel performance code: Normal operating conditions (part 1),'' {\em
  Nuclear Engineering and Design}, vol.~282, pp.~158--168, 2015.

\bibitem{Prudil15b}
A.~Prudil, B.~Lewis, P.~Chan, J.~Baschuk, and D.~Wowk, ``Development and
  testing of the {FAST} fuel performance code: Transient conditions (part 2),''
  {\em Nuclear Engineering and Design}, vol.~282, pp.~169--177, 2015.

\bibitem{Long02}
Y.~Long, Y.~Yuan, M.~Kazimi, R.~Ballinger, and E.~Pilat, ``A fission gas
  release model for high-burnup {LWR} {ThO$_2$-UO$_2$} fuel,'' {\em Nuclear
  Technology}, vol.~138, no.~3, pp.~260--272, 2002.

\bibitem{Lee04}
C.~Lee, Y.~Yang, Y.~Kim, D.~Kim, and Y.~Jung, ``Irradiation behavior of
  thoria-urania fuel in a {PWR},'' {\em Nuclear Technology}, vol.~147, no.~1,
  pp.~140--148, 2004.

\bibitem{Bjork15}
K.~Bj{\"{o}}rk and L.~Kekkonen, ``Thermal-mechanical performance modeling of
  thorium-plutonium oxide fuel and comparison with on-line irradiation data,''
  {\em Journal of Nuclear Materials}, vol.~467, pp.~876--885, 2015.

\bibitem{Boer16}
B.~Boer, S.~Lemehov, M.~Weber, Y.~Parthoens, M.~Gysemans, J.~McGinley,
  J.~Somers, and M.~Verwerft, ``Irradiation performance of {(Th,Pu)O$_2$} fuel
  under pressurized water reactor conditions,'' {\em Journal of Nuclear
  Materials}, vol.~471, pp.~97--109, 2016.

\bibitem{Booth57}
A.~Booth, ``A method of calculating fission gas release from uo$_2$ fuel and
  its application to the {X-2-f} loop test,'' Tech. Rep. 469, Atomic Energy of
  Canada Limited, 1957.

\bibitem{White83}
R.~White and M.~Tucker, ``A new fission-gas release model,'' {\em Journal of
  Nuclear Materials}, vol.~118, no.~1, pp.~1--38, 1983.

\bibitem{Karam08}
M.~Karam, F.~Dimayauga, and J.~Montin, ``Fission gas release of {(Th, Pu)O$_2$
  CANDU} fuel,'' in {\em Proceedings of the 10th International Conference on
  CANDU Fuel}, (Ottawa, ON, Canada), 2008.

\bibitem{Karam10}
M.~Karam, F.~Dimayauga, and J.~Montin, ``Post-irradiation examination of candu
  fuel bundles fuelled with {(Th,Pu)O$_2$},'' in {\em Proceedings of the 11th
  International Conference on {CANDU} Fuel}, (Niagara Falls, ON, Canada), 2010.

\bibitem{Floyd11}
M.~Floyd, ``Advanced fuel cycle development at {Chalk River Laboratories},'' in
  {\em Proceedings of the International Conference on the Future of Heavy Water
  Reactors}, (Ottawa, ON, Canada), 2011.

\bibitem{Floyd16}
M.~Floyd, R.~Beier, J.~Bell, A.~Bergeron, E.~Bulemela, G.~Cota-Sanchez, R.~S.
  Dickson, K.~Leeder, S.~Livingstone, A.~Quastel, and M.~Saoudi, ``Progress in
  the development of thoria fuel science \& technology,'' in {\em 13th
  International Conference on CANDU Fuel}, (Kingston, ON, Canada), 2016.

\bibitem{Long04}
Y.~Long, L.~Siefken, P.~Hejzlar, E.~Loewen, J.~Hohorst, P.~MacDonald, and
  M.~Kazimi, ``The behavior of {ThO$_2$}-based fuel rods during normal
  operation and transient events in {LWR}s,'' {\em Nuclear Technology},
  vol.~147, no.~1, pp.~120--139, 2004.

\bibitem{Sherman50}
J.~Sherman and W.~Morrison, ``Adjustment of an inverse matrix corresponding to
  a change in one element of a given matrix,'' {\em Annals of Mathematical
  Statistics}, vol.~21, no.~1, pp.~124--127, 1950.

\bibitem{Transtrum12}
M.~Transtrum and J.~Sethna, ``Improvements to the {L}evenberg-{M}arquardt
  algorithm for nonlinear least-squares minimization,'' {\em eprint
  arXiv:1201.5885}, vol.~http://arxiv.org/abs/1201.5885, 2012.

\bibitem{Piro13}
M.~Piro, S.~Simunovic, T.~Besmann, B.~Lewis, and W.~Thompson, ``The
  thermochemistry library \textsc{thermochimica},'' {\em Computational
  Materials Science}, vol.~67, pp.~266--272, 2013.

\end{thebibliography}

\end{document}